\documentclass{elsart}
\usepackage{amssymb,amsmath,graphicx}

\begin{document}
\begin{frontmatter}

\title{Dynamics near nonhyperbolic fixed points or nontransverse homoclinic points}

\journal{Mathematics and Computers in Simulation}

\author[spbu]{Sergey Kryzhevich},
\corauth[cor]{Corresponding author. Tel.:+7 921 918 10 76; fax: +7 812 428 69 90.}
\ead{kryzhevitz@rambler.ru}
\author[spbu]{Sergei Pilyugin}
\ead{sergeipil47@mail.ru}

\address[spbu]{Faculty of Mathematics and Mechanics, Saint-Petersburg State University, 28, Universitetsky pr., Peterhof, Saint-Petersburg, Russia, 198503}

\begin{abstract}

We study dynamics in a neighborhood of a nonhyperbolic fixed point or an irreducible homoclinic tangent point. General type conditions for the existence of infinite sets of periodic points are obtained. A new method, based on the study of the dynamics of center disks, is introduced. Some results on shadowing near a non-hyperbolic fixed point of a homeomorphism are obtained.
\end{abstract}

\begin{keyword}
Partial hyperbolicity \sep center unstable manifold \sep homoclinic point \sep shadowing
\MSC 37B10 \sep 37B25 \sep 37G30
\end{keyword}
\end{frontmatter}

\section{Introduction}

Many important problems of the bifurcation theory or the theory of strongly nonlinear and discontinuous dynamical systems can be reduced to the problem on topological structure of trajectories in a neighborhood of a nonhyperbolic fixed point. Started from works of Lyapunov and Poincar\'e, the theory of non-hyperbolic systems had a great breakthrough due to the appearance of the reduction principle \cite{krpil1} (see also \cite{krpil2,krpil3,krpil4,krpil5}).

Nonlinear phenomena, for example, chaotic dynamics, are possible in a neighborhood of a nonhyperbolic fixed point \cite{krpil6,krpil7,krpil8,krpil9,krpil10,krpil12,krpil11,krpil14,krpil13,krpil15,krpil16,krpil17,krpil18} (see also references therein and a closely related result of \cite{krpil19}). Bo Deng \cite{krpil8} has described the case of a transverse homoclinic point corresponding to a unique eigenvalue on the unit circle (i.e., the center manifold is one-dimensional). Making additional assumptions on the smoothness of the diffeomorphism and on its Jacobi matrix at zero, he established the existence of analogues of the Smale horseshoe in a neighborhood of a homoclinic point. Moreover, the appearing transitive invariant sets persist while parameters of the system are slightly perturbed, though the local structure of the set of nonwandering points may change sufficiently [20].

Basing on an example of R.\, Ma\~ne \cite{krpil15}, Buzzi and Fisher \cite{krpil9} have described some properties of transitive non-Anosov diffeomorphisms of smooth manifolds. Sufficient results on existence of nonhyperbolic chaos in Hamiltonian systems have been obtained.

The structure of invariant sets for mappings, represented as skew products over the Smale horseshoe has been discussed in the paper [11] and successive works. It was shown that there is an open set of diffeomorphisms of that type, having infinite sets of sources and sinks.

A review of other recent results can be found in \cite{krpil16}.

An important result involving properties of so-called weakly hyperbolic invariant sets (particularly, their structural stability) has been obtained in the paper \cite{krpil17}.

The appearing transitive invariant sets are not hyperbolic, so their persistence is, at least, not evident. To study this problem, the theory of normal and partial hyperbolicity \cite{krpil6,krpil21,krpil22,krpil23,krpil26,krpil25} may be applied.

The dynamics in a neighborhood of a nonhyperbolic fixed point is related to one in the neighborhood of a homoclinic tangency \cite{krpil27,krpil28,krpil29,krpil30,krpil31,krpil33,krpil32,krpil11,krpil34,krpil35,krpil36,krpil37} (see also references therein). Usually, we can reduce one problem to another by using a transformation of variables, that is not smooth at points of the stable manifold.

One of the most interesting results in this area has been obtained by Newhouse \cite{krpil35} and later on developed in the papers by Gonchenko, Shil'nikov, Turaev, Vasil'eva, and others \cite{krpil29,krpil30,krpil31,krpil33,krpil32,krpil34,krpil35,krpil36,krpil37}. A neighborhood of a homoclinic tangency may contain an infinite set of periodic attracting points. The Lyapunov exponents of these periodic points tend to $0$.

However, there is no general description of the dynamics in a neighborhood of a nonhyperbolic fixed point.

The rest of the article consists of five sections where the new results are presented, conclusions and appendix. In Section 2, we establish a criterium on the existence of infinitely many periodic points in a neighborhood of a nonhyperbolic fixed point. We suppose, that there exists a Lyapunov function that does not allow leaves of the center unstable manifold to shrink. For this case, one can prove an analogue of the $\lambda$~-- lemma (see the appendix) and the uniqueness of the center unstable manifold and generalize the Smale - Birkhoff theorem. A different approach to application of Lyapunov functions for the nonhyperbolic case was developed in \cite{krpil38,krpil39}.

In Section 3, we apply the obtained results to study a special case of homoclinic tangency (for example, a cubic type tangency for 2D systems). In this case, we obtain an infinite set of periodic points, as well.

In Section 4, we use the Lyapunov function approach to prove a new criterium of shadowing in a neighborhood of a nonhyperbolic fixed point. The proof of our shadowing result uses ideas similar to those applied in \cite{krpil43}. The advantage of our approach is that we work with homeomorphisms and do not refer to smoothness of the dynamical system.

In Section 5, an example illustrating the results of the previous section is given.

\section{Nonhyperbolic homoclinic points}

Let $Q$ be a bounded domain in the Euclidean space ${\mathbb R}^n$ containing the origin. Consider the set $X$ of $C^1$ diffeomorphisms
$F:{\mathbb R}^n\to {\mathbb R}^n$ such that $F(\overline Q)\subset Q$, where $\overline Q$ is the closure of the set $Q$. We identify all the diffeomorphisms that coincide on $\overline Q$.

We endow $X$ with the $C^1$~-- metrics defined by the formula
$$d_1(F,F')=\max_{x\in \overline{Q}}|F(x)-F'(x)|+\max_{x\in \overline{Q}}|DF(x)-DF'(x)|.$$
Here and below, we denote by the symbol $|\cdot|$ both the Euclidian norm of a vector and the corresponding operator matrix norm and by $DF$ and $DF'$ the Jacobi matrices of the corresponding mappings. Fix a diffeomorphism $F\in X$ and assume that the point $0$ is a fixed point of $F$. Let $\chi_1,\ldots,\chi_n$ be the eigenvalues of the matrix $A=DF(0)$ (some of them may be equal). Assume that
$$|\chi_1|\leq |\chi_2|\leq\ldots\leq |\chi_{s}|<1\leq |\chi_{s+1}|\leq\ldots\leq |\chi_n|$$
for some $s\in \{1,2,\ldots, n\}$.

Let $u=n-s$. Consider the eigenspace ${\cal S}$ corresponding to $\chi_1$, \ldots, $\chi_{s}$ and the eigenspace ${\cal U}$ corresponding to the eigenvalues $\chi_{s+1}$, \ldots, $\chi_n$. Following the paper \cite{krpil3}, we call the spaces ${\cal S}$ and ${\cal U}$ the stable and center unstable spaces, respectively.

Without loss of generality, we assume that the space $\cal S$ is the linear hull of the first $s$ coordinate vectors and $\cal U$ is the linear hull of the last $u$ ones.

Here we quote some definitions and results from the book \cite{krpil22}, changing some notation. Let $\rho\in (0,1]$.

\noindent\textbf{Definition 2.1} \emph{\cite[p.58]{krpil22}. A linear endomorphism of a Banach space $L:E\to E$ is $\rho$-pseudohyperbolic if its spectrum lies off of the circle of radius $\rho$.}

A splitting $$E=E^s\oplus E^{cu}$$ is called \emph{canonical} for a $\rho$-pseudohyperbolic endomorphism $L:E\to E$ if both subspaces $E^s$ and $E^{cu}$ are $L$~-- invariant, the spectrum of $L_{s}=L|_{E^{s}}$ lies inside the open ball of radius $\rho$ and the spectrum of $L_{cu}=L|_{E^{cu}}$ lies off this ball.

\noindent\textbf{Theorem 2.1} \emph{\cite[p.58]{krpil22}. Let $L:E\to E$ be a $\rho$-pseudohyperbolic endomorphism of a Banach space, $E=E^{s}\oplus E^{cu}$ be the canonical splitting, $F:E\to E$ be a $C^1$ map, $F(0)=0$,
$$\begin{array}{c} F(x)=Lx+f(x), \quad \mbox{and} \quad |f(x_1)-f(x_2)|\le\varepsilon|x_1-x_2|, \quad x_{1,2}\in E. \end{array} \eqno (2.1)$$
Then there exists an $\varepsilon_0$ such that if $\varepsilon<\varepsilon_0$, then the sets $W^{s}$ and $W^{cu}$, called the stable and center unstable manifolds, respectively, and defined by
$$\begin{array}{l}
W^s=\bigcap_{m=0}^\infty F^{-m}(C_1), \qquad C_1=\{(y,z)\in E^s\oplus E^{cu}:|y|\ge |z|\};\\
W^{cu}=\bigcap_{m=0}^\infty F^m(C_2), \qquad C_2=\{(y,z)\in E^s\oplus E^{cu}:|y|\le |z|\};
\end{array}$$
are graphs of $C^1$ maps $E^s\to E^{cu}$ and $E^{cu}\to E^s$. They are characterized by
$$\begin{array}{l}
x\in W^s \Leftrightarrow |F^m(x)|\rho^{-m} \to 0 \quad \mbox{as}\quad m\to +\infty; \\[5pt]
x\in W^{cu} \Leftrightarrow \mbox{there exist all inverse images}\quad F^{-m}(x), \\
|F^{-m}(x)|\rho^{-m} \to 0 \quad \mbox{as}\quad m\to +\infty,\end{array}$$
and "$\to 0$"\ may be replaced by "stays bounded". If $\|L_s\|\|L_u^{-1}\|<1$
then both manifolds $W^s$ and $W^{cu}$ are $C^1$. The manifolds $W^s$ and $W^{cu}$ continuously depend on $F$ in the $C^1$ sense.}

The mapping $F$ may not satisfy the estimates (2.1) globally (actually, we even do not assume that this mapping is globally defined). In this case the locally invariant manifolds $W^{s}_{loc}$ and $W^{cu}_{loc}$ can still be defined, but $W_{loc}^{cs}$ is, in general, not unique any more, see [3] and [22, \S 5]. The sets $W^s_{loc}$ and $W^{cu}_{loc}$ are usually called the local stable and local center unstable manifolds of the fixed point 0.

If $Dg(0)=0$, the manifold $W^s$ is tangent to $E^s$ and any of the manifolds $W^{cu}$ is tangent to $E^{cu}$ at 0.

Similarly, one may define $\rho$-pseudohyperbolicity for $\rho>1$. For this case, an analogue of Theorem 2.1 is true and there are two invariant manifolds $W^{cs}$ and $W^u$, called the center stable and unstable manifolds, respectively.

If the conditions of Theorem 2.1 are satisfied locally and, consequently, there exist local manifolds $W^s_{loc}$ and $W^{cu}_{loc}$, we can extend these manifolds to invariant sets
$$W^{\sigma}=\bigcup_{k\in {\mathbb Z}} F^k(W^\sigma_{loc}), \qquad \sigma\in \{s,cu\},$$
still called the stable and center unstable manifold of 0, respectively.

The next definition is a new one.

\noindent\textbf{Definition 2.2.} \emph{(Fig.\,1). We say that the fixed point $0$ of the mapping $F$ is \emph{strongly conditionally unstable} if we can choose a manifold $W^{cu}_{loc}$ so that
\begin{itemize}
\item[1)] there exists a neighborhood $U_0$ of the origin and a continuous mapping $V:U_0\to [0,+\infty)$
such that $V(x)=0$ if and only if $x\in W^s_{loc}$, and $V(F(x))\geq V(x)$
for all $x\in U_0\bigcap F^{-1}(U_0)$;
\item[2)] the fixed point $0$ of the restriction $F^{-1}|_{W^{cu}_{loc}}$ is asymptotically stable.
\end{itemize}}

\begin{figure}\begin{center}
\includegraphics*[width=2in]{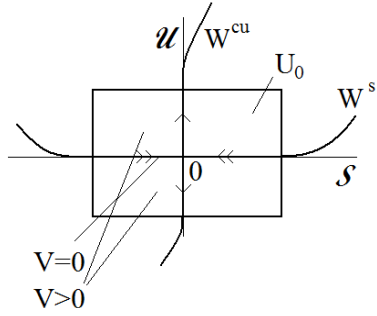}
\end{center}
\caption{\it The stable and center unstable manifolds.}
\end{figure}

\noindent\textbf{Remark 2.1.} \emph{It follows from \cite[Theorem 5A.3]{krpil22}, that the strong conditional instability implies the uniqueness of the center unstable manifold.}

We say that the fixed point $0$ of the mapping $F$ is \emph{strongly conditionally stable} if $0$ is strongly conditionally unstable for the diffeomorphism $F^{-1}$. In this case,
$$W^{cs}=W^{cs}(F)=W^{cu}(F^{-1}), \quad \mbox{and} \quad W^{u}=W^{u}(F)=W^{s}(F^{-1}).$$

Note some obvious properties of strong conditional (in)stability.
\begin{enumerate}
\item If the diffeomorphism $F(x)$ can be locally represented as $F=F_s\times F_u$, where the point $0\in {\mathbb R}^{s}$ is a hyperbolic attracting fixed point for the diffeomorphism $F_s:{\cal S}\to {\cal S}$ and the point $0\in {\cal U}$ is a repeller for the diffeomorphism $F_u:{\cal U}\to {\cal U}$, then the fixed point 0 of the mapping $F$ is strongly conditionally unstable.
\item Strong conditional (in)stability is invariant with respect to homeomorphic transformations of coordinates.
\item Every hyperbolic fixed point is strongly conditionally stable and strongly conditionally unstable at the same time.
\end{enumerate}

In what follows, we assume that one of the following symmetric conditions is satisfied.

\noindent\textbf{Condition 2.1.} \emph{The point $0$ is a strongly conditionally stable fixed point of the diffeomorphism $F$. There exist disks $w^{cs}\subset W^{cs}$ and $w^{u}\subset W^u$ that intersect transversally at a point $p\neq 0$.}

\noindent\textbf{Condition 2.2.} \emph{(Fig.\,2.) The point $0$ is a strongly conditionally unstable fixed point of the diffeomorphism $F$. There exist disks $w^{cu}\subset W^{cu}$ and $w^{s}\subset W^s$, that intersect transversally at a point $p\neq 0$.}

\begin{figure}\begin{center}
\includegraphics*[width=2in]{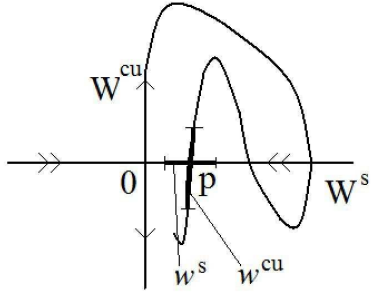}
\end{center}
\caption{\it Non-hyperbolic homoclinic intersection.}
\end{figure}

\noindent\textbf{Theorem 2.4.} \emph{Let $F\in X$ and let either Condition 2.1 or Condition 2.2 be satisfied. Then for any neighborhood $U$ of the origin there exists a number $\delta>0$ such that for any $G\in X$ with
$$d_1(F,G)<\delta \eqno (2.2)$$
there is an infinite subset $P_G\subset U$ with the following properties.
\begin{itemize}
\item[1)] Every point $q\in P_G$ is a periodic point of $G$;
\item[2)] for any $m\in {\mathbb N}$ there is a point $q\in P_G$ such that the minimal period of $q$ is larger than $m$;
\item[3)] ${\mathop{\rm card}\ } \overline P_G=\aleph$.
\end{itemize}}

\noindent\textbf{Remark 2.2.} \emph{The sets $P_G$ corresponding to different mappings $G$ may have different topological structures.}

\textbf{Proof.} In the proof, we consider the case where Condition 2.2 is satisfied.

Without loss of generality, we may assume that $A=\mathop{\rm diag} (B,C)$, where the $s\times s$ matrix $B$ and the $u\times u$ matrix $C$ are such that $|B|=a_0<1$ and $b_0=1/|C^{-1}|>a_0$ (we can satisfy the second inequality due to Condition 2.2). We represent $x=(y,z)$, where the vector $y$ consists of the first $s$ components of the vector $x$. Then ${\cal S}=\{(y,0): y\in {\mathbb R}^s\}$ and ${\cal U}=\{(0,z):z\in {\mathbb R}^{u}\}$. Moreover, we may assume that $W^s_{loc}\subset {\cal S}$ and $W^{cu}_{loc} \subset {\cal U}$.

Since the trajectory of the point $p$ consists of homoclinic points with similar properties, we may assume that $p\in W^s_{loc}$ and $p=(y_p,0)$.
Let $U^0$ be a neighborhood of the origin, so small that $p\notin \overline{U^0}$,
$$W^s_{loc}\bigcap U^0=\{(y,0)\in U^0\} \quad \mbox{and}\quad W^{cu}_{loc}\bigcap U^0=\{(0,z)\in U^0\}.$$
We take standard coordinates $(y,z)$ in the neighborhood $U^0$. Since $w^s$ and $w^{cu}$ are transverse at the point $p$, we can introduce coordinates $(y,z')$ in a neighborhood $U^1$ of $p$ such that $U^0\cap U^1=\emptyset$, the transformation $(y,z)\leftrightarrow (y,z')$ is smooth and a small disk in $w^{cu}$ containing $p$ is given in coordinates $(y,z')$ by the equality $y=y_p$.

Consider two disjoint neighborhoods $U_0$ and $U_1$ of the points $0$ and $p$ respectively, given by
$$\begin{array}{c} U_0=\{x=(y,z):|y|\leq \varepsilon^y_0,\quad |z|\leq \varepsilon^z_0\};\\
U_1=\{x=(y,z):|y-y_p|\leq \varepsilon^y_1,\quad |z'|\leq \varepsilon^z_1\}.\end{array}$$
Here $\varepsilon^y_{0,1}$ and $\varepsilon^z_{0,1}$ are positive and so small that $U_i\subset U^i$, $i=0,1$.

We define an \emph{admissible disk at $U_0$} as a set $D=\{(\eta(z),z):|z|\leq \varepsilon^z_0\}$.
Here the function
$\eta\in C^1(B^{cu}\to {B^s})$ is $C^1$--smooth,
$$\max_{|z|\leq \varepsilon_0^z}|\eta(z)|\leq \varepsilon_0^y \quad \mbox{and} \quad
\max_{|z|\leq \varepsilon_0^z}|D\eta(z)|\leq 1.$$

Similarly, we define an \emph{admissible disk at $U_1$} as a set
$D=\{(\eta(z'),z'):|z'|\leq \varepsilon^z_1\}$.
where the function $\eta$ is $C^1$--smooth,
$$\max_{|z'|\leq \varepsilon_1^z}|\eta(z')-y_p|\leq \varepsilon_1^y \quad \mbox{and} \quad \max_{|z'|\leq \varepsilon_1^z}|D\eta(z')|\leq 1.$$
Let ${\cal D}_i$ be sets of disks admissible at $U_i$, $i=0,1$.
We denote ${\cal D}={\cal D}_0\bigcup {\cal D}_1$ and endow both sets ${\cal D}_i$ with the natural $C^1$ metrics ${\mathop {\rm dist}}_1$, identifying admissible disks with corresponding maps $\eta$.

We use the following statement, which is very similar to the well known $\lambda$~- lemma, see another generalization in \cite{krpil40}.

\noindent\textbf{Lemma 2.1.} \emph{Assume that $0$ is a strongly conditionally unstable fixed point of the mapping $F$. Then there exists a smooth disk $N\subset W^{cu}_{loc}$ containing the point $0$ and such that for every positive $\varepsilon$ there exists $m(\varepsilon)\in {\mathbb N}$ having the following property. For every $m\geq m(\varepsilon)$, $D\in {\cal D}$ there exists an embedding $h_m$ of the disk $N$ into ${\mathbb R}^n$ such that $h_m(N)\subset F^m(D)$ and $\mbox{\rm dist}_1(h_m,{\mathop{\rm id}})<\varepsilon$.}

A proof of this statement is given in the appendix.

We can take $k$ so large and $\delta$ so small that the following statement holds.
\begin{enumerate}
\item For any $G$ satisfying (2.2) and for $i,j\in\{0,1\}$, the intersections $$G^k(U_i)\bigcap U_j$$ contain components $U_{ij}$ such that if
$D\in {\cal D}_i$, then the intersection $$S_j(D)=G^k(D)\bigcap U_{ij}$$ is an admissible disk at $U_j$.
\item $\mbox{\rm dist}_1(S_j(D),S_j(D'))\leq \mbox{\rm dist}_1(D,D')/2$ for $i,j\in \{0,1\}$ and all $$D,D'\in {\cal D}_i.$$
\end{enumerate}

This can be done due to Lemma 2.1. We apply this lemma to the mapping $F$ and then note that if the components $U_{ij}$ exist for $F$, they can be chosen for all $G$ sufficiently $C^1$--close to $F$ since $k$ is fixed.

Then for any admissible disk $D\subset U_i$, the admissible disk $S_j(D)\subset G^k(D)\bigcap U_{ij}$
is uniquely defined and depends continuously on $D$ in the metrics $\mathop{\rm dist}_1$.

Consider the set $\Sigma$, consisting of infinite one-side sequences $$a=\{a_{k}\in\{0,1\}:k\in {\mathbb Z}^+\}.$$
Let us define a metrics in the set $\Sigma$ in the standard way:
$$d(a,b)=\sum_{k=0}^\infty 2^{-k}|a_k-b_k|.$$
We identify periodic subsequences of $\Sigma$, that can be obtained by an infinite repetition of finite sequences $(a_0,\ldots, a_N)$,
$a_j\in \{0,1\}$, with these finite sequences. If a finite sequence may be obtained from another one by finite repetition, we also identify these sequences.

To any periodic sequence $a$ generated by a set $\{a_0,\ldots, a_N\}$ we assign the admissible disk $D_a$ that is the unique fixed point of the contracting mapping $S_{a_0}\circ S_{a_1}\circ \ldots\circ S_{a_N}$.
Then the following inclusions hold:
$$\begin{array}{c}D\in U_{a_0}, \quad D\in S_{a_0}({\cal D}_{a_1}), \quad D\in S_{a_0}\circ S_{a_1} ({\cal D}_{a_2}),\ldots,\\ D \in S_{a_0}\circ S_{a_1}\circ\ldots\circ S_{a_{N-1}} ({\cal D}_{a_N}).\end{array}\eqno (2.3)$$
Properties of the mappings $S_i$ imply that there exists a constant $C>0$ such that if $k\in {\mathbb N}$ and the first $k$ entries of finite sequences
$a=(a_0,a_1,\ldots,a_{N_1})$ and $b=(b_0,b_1,\ldots,b_{N_2})$ coincide, then
$$\mbox{\rm dist}_1 (D_a,D_b)\le C2^{-k}.$$

Consequently, for any converging sequence $\{b^k\}\subset \Sigma$ consisting of periodic elements, the corresponding sequence $D_{b^k}$ converges in the space $\cal D$. Fix an element $a=(a_0,a_1,\ldots, a_N,\ldots)\in \Sigma$
and denote $a^k=(a_0,\ldots, a_k)$. Clearly, $a^k\to a$. Let $D_a=\lim D_{a^k}$. Consider the set $K=\{D_a: a\in \Sigma\}$.
Let $H:\Sigma \to {\cal D}$ be the parameterizing mapping: $H(a)=D_a$. Denote by $\sigma_i$ the adding of $i\in\{0,1\}$ to the left-hand side of an element $a\in\Sigma$. For an arbitrary $a\in \Sigma$ we denote $1a=\sigma_1(a)$ and $0a=\sigma_0(a)$.

Let us prove that the disks corresponding to different elements of the space $\Sigma$ are distinct. Let $a,b\in \Sigma$, $a\neq b$. Consider the least integer $j$ such that $a_j\neq b_j$. If $j=0$, the disks $D_a$ and $D_b$ appertain to distinct sets $U_i$ since all the approximating disks do. Otherwise, due to (2.3), there exist the inclusions
$$\begin{array}{c}D_a \in S_{a_0}\circ S_{a_1}\circ\ldots\circ S_{a_{j-1}} ({\cal D}_{a_j}), \\ D_b \in S_{b_0}\circ S_{b_1}\circ\ldots\circ S_{b_{j-1}} ({\cal D}_{b_j})= S_{a_0}\circ S_{a_1}\circ\ldots\circ S_{a_{j-1}} ({\cal D}_{b_j}).\end{array}$$
Hence the disks $D_a$ and $D_b$ do not intersect.

\noindent\textbf{Lemma 2.2.} \emph{For $i\in \{0,1\}$, $\sigma_i\circ H=H\circ S_i$.}

\textbf{Proof.} Let $i=0$, the case $i=1$ is similar. Fix a sequence $a=\{a_k:k\in{\mathbb Z}^+\}\in \Sigma$ and the corresponding disk
$$\begin{array}{c} D_a= {\cal D}_{a_0}\bigcap {\cal D}_{a_0}(U_{a_1})\bigcap {\cal D}_{a_0}\circ S_{a_1} ({\cal D}_{a_2})\bigcap \ldots \bigcap \\ S_{a_0}\circ S_{a_1}\circ\ldots\circ S_{a_{N-1}} ({\cal D}_0{a_N})\bigcap \ldots\end{array}$$

Then
$$\begin{array}{c} D_{0a}= {\cal D}_0 \bigcap S_0 ({\cal D}_{a_0})\bigcap S_0\circ S_{a_0}({\cal D}_{a_1})\bigcap \ldots \\
\bigcap S_0\circ S_{a_0}\circ S_{a_1}\circ\ldots\circ S_{a_{N-1}} ({\cal D}_{a_N})\bigcap \ldots=\\
S_0 ({\cal D}_{a_0})\bigcap S_0\circ S_{a_0}({\cal D}_{a_1})\bigcap \ldots \\
\bigcap S_0\circ S_{a_0}\circ S_{a_1}\circ\ldots\circ S_{a_{N-1}} ({\cal D}_{a_N})\bigcap \ldots=
S_0 (D_a).\quad \blacksquare\end{array}$$

Periodic points of the shift mapping $\sigma$ are dense in $\Sigma$, and there is a point
$$a^*=\{a_k^*,k\in {\mathbb Z}^+\}\in \Sigma$$
whose positive semi-orbit is dense \cite{krpil41}. We call an admissible disk periodic if it corresponds to a periodic sequence. These periodic admissible disks are dense in $K$ and there is a dense sequence $\{D_k:k\in {\mathbb Z}^+\}\subset K$ such that
$$S_{a_k^*}(D_k)=D_{k-1}$$
for all $k\in {\mathbb Z}^+$.

Note that for any periodic disk $D$ of period $m$, $D\subset G^{km}(D)$.
Then due to the Brauer's theorem, there is a point $x\in D$ such that $x=G^{km}(x)$.

Consider a convergent sequence $\{a_k\}\subset\Sigma$ such that every $a_k$ is periodic. Denote $a=\lim a_k$. Let $x_k\in D_{a_k}$ be a sequence of periodic points. There exists a limit point $x\in D_a$ of the sequence $x_k$. Obviously, $x\in \overline{P_G}$. Hence, $D_a\bigcap \overline{P_G}\neq \emptyset$ for any $a\in\Sigma$. This proves that ${\mathop{\rm card}\ } \overline P_G=\aleph$. $\blacksquare$

\section{Nontransverse homoclinic points}

We say that two $C^1$~-- smooth submanifolds $W^s$ and $W^u$ of the Euclidean space ${\mathbb R}^n$ intersect \emph{quasitransversally} at a point $p$ if
$\dim W^s+\dim W^u=n$ and there exist a neighborhood $U$ of the point $p$ and a $C^1$ smooth coordinate system $\xi={\mathop{\rm col}}(\eta,\zeta)$ in $U$ with the following properties (Fig.\, 3).

\begin{figure}\begin{center}
\includegraphics*[width=2in]{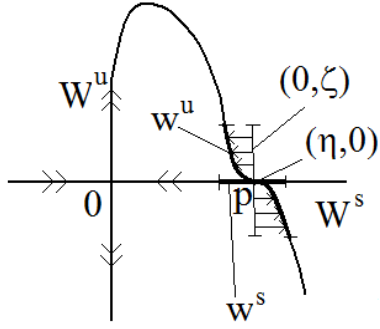}
\end{center}
\caption{\it Quasitransverse homoclinic tangency.}
\end{figure}

Denote by $w^s$ and $w^u$ the connected components of intersections of $W^s\bigcap U$ and $W^u\bigcap U$, containing the point $p$.
\begin{enumerate}
\item $\dim \zeta=\dim W^u$;
\item the manifolds $\{(0,\zeta)\}$ and $w^s$ intersect transversally at the point $p$;
\item the mapping $\xi:U\to \xi(U)$ is a local diffeomorphism;
\item $\zeta(x)=0$ for any $x\in w^s$;
\item there exists a number $\delta>0$ such that the set $w^u$ is the graph of a function $\eta=g(\zeta)$, $|\zeta|<\delta$ such that the mapping $g$ is smooth for all $\zeta: 0<|\zeta|<\delta$.
\end{enumerate}

We say that two subsets $X^1$ and $X^2$ of the Euclidean space ${\mathbb R}^n$ intersect quasitransversally at a point $p$ if there exists a neighborhood $U$ of the point $p$ such that the connected components $W^1$ and $W^2$ of intersections $X^1\bigcap U$ and $X^2\bigcap U$ that contain the point $p$ are $C^1$ smooth disks and intersect quasitransversally at this point.

\noindent\textbf{Remark 3.1.} \emph{As we show later, this condition implies that there is a locally smooth invertible (but not always diffeomorphic) transformation of coordinates, which makes a quasitransverse intersection transverse.}

For example a cubic type tangency of two curves in ${\mathbb R}^2$ is a quasitransverse intersection.

\noindent\textbf{Theorem 3.1.}\emph{ Let $F\in X$ be such that $x=0$ is a hyperbolic fixed point. Assume that $F$ can be $C^1$ linearized in a neighborhood $U$ of 0. Assume that the corresponding stable and unstable manifolds ($W^s$ and $W^u$) intersect quasitransversally at a point $p\neq 0$.
Then for any neighborhood ${\cal V}$ of the origin there is an infinite subset $\Pi\in V$ with the following properties:
\begin{enumerate}
\item[1)] every point $q\in \Pi$ is a periodic point of $F$;
\item[2)] for any $m\in {\mathbb N}$ there is a point $q\in \Pi$ such that the minimal period of $q$ is larger than $m$;
\item[3)] ${\mathop{\rm card}\ } \overline \Pi=\aleph$.
\end{enumerate}}

\textbf{Proof.} We assume without loss of generality that there is a small neighborhood ${\cal V}$ of the origin such that $F|_V$ is a linear mapping of the form $F(x)=Ax$, where $A={\mathop {\rm diag}}(B,C)$, $B$ is an $s\times s$ matrix, $C$ is a $u\times u$ matrix ($s+u=n$), $|B|<1$ and $|C^{-1}|<1$.

Also, we can assume the following:
\begin{enumerate}
\item The neighborhood $U$ in the definition of quasitransversality of intersection of $W^s$ and $W^u$ at the point $p$ is a subset of ${\cal V}$;
\item $p=(y_p,0)$;
\item the local invariant manifolds $W^s_{loc}$ and $W^u_{loc}$ defined in a neighborhood of the origin are given by $z=0$ and $y=0$ respectively,
\item $w^s$ is given by $z=0$,
\item $w^u$ is given by $y=y_p$ and
\item $F(y,z)=(F_1(y),F_2(z))=(By, Cz)$ for all $x=(y,z)\in U$.
\end{enumerate}

Considering, if necessary, the mapping $F^2$ instead of $F$ and the matrix $C^2$ instead of $C$, we may assume that there exists a real-valued matrix $P$ such that $C=\exp(P)$. All the eigenvalues of the matrix $P$ have positive real parts. Making, if necessary, a linear transformation of variables, we may suppose, that the unit sphere $\{z\in {\mathbb R}^n:|z|=1\}$ is transverse to trajectories of the system
$$\dot z=Pz. \eqno (3.1)$$

The main idea of the proof is to construct a transformation $(y,z)\to (y,\hat z)$ in order to obtain a transverse nonhyperbolic homoclinic point instead of the quasitransverse tangency. Then the result of Theorem 2.4 can be applied.

We start with the following obvious technical statement.

\noindent\textbf{Lemma 3.2.}\emph{ Let $\delta:[0,\rho)\to [0,+\infty)$ be a continuous nondecreasing function such that $\delta(0)=0$ and $\delta(\xi)>0$ for all $\xi>0$. Then there exists a $C^{\infty}$ smooth function $\delta_0:[0,\rho)\to [0,+\infty)$ such that $\delta_0(0)=0$, $\delta_0'(\xi)\geq 0$, $0<\delta_0(x)\leq \delta(x)$ for all $\xi>0$, and $\delta_0^{(k)}(0)=0,\qquad  k\in {\mathbb N}$.}



Figure 4 illustrates a possible way of constructing the function $\delta_0$.

\begin{figure}\begin{center}
\includegraphics*[width=2in]{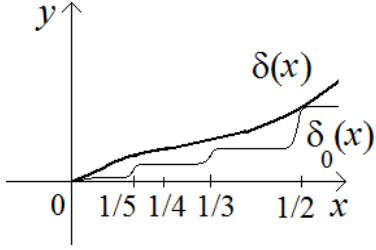}
\end{center}
\caption{\it Functions $\delta$ and $\delta_0$.}
\end{figure}

\noindent\textbf{Lemma 3.3. }\emph{There exists a positive number $\rho$ and a $C^\infty$ smooth function $t:(0,\rho)\to [0,+\infty)$ with the following properties.
\begin{enumerate}
\item $t'(\xi)<0$ for all $\xi\in (0,\rho)$.
\item $\lim\limits_{\xi\to 0} \xi^2 t(\xi)=+\infty$ and $\lim\limits_{\xi\to 0} t'(\xi)=-\infty$.
\item Let a mapping $h$ be defined by the formula
$$h(\hat z)=\exp(-t(|\hat z|^2)P)\hat z \eqno (3.2)$$
for $\hat z\neq 0$ and $h(0)=0$. Then $h$ is a homeomorphism of the domain $U=\{\hat z:|\hat z|\leq \rho\}$ to $h(U)$.
\item For any $\hat z$ such that $0<|\hat z|\leq \rho$, the Jacobi matrix $Dh(\hat z)$ is well-defined and invertible.
\item If $\hat g(\hat z)=g(h(\hat z))$, then $D\hat g(0)=0$.
\item Let $\hat F_2=h^{-1}\circ F_2 \circ h$. Then the mapping $\hat F_2$ is smooth in a neighborhood of zero, and $D\hat F_2 (0)=E$.
\end{enumerate}}

\textbf{Proof}. Let $\delta$ be a continuous monotonous function such that $\delta(\varepsilon)>0$ for any $\varepsilon>0$ and the inequality $|z|\leq \delta(\varepsilon)$ provides that $|g(z)|\leq \varepsilon$. Select a $C^\infty$ smooth function $\delta_0(\varepsilon)$ according to Lemma 3.2. Let positive constants $\chi$ and $K$ be such that
$$|\exp(-Pt)|\leq K\exp (-\chi t)\eqno (3.3)$$
for any $t>0$.

Let $$t(\xi)=\dfrac1{\delta_0(\xi)}+\exp\left(\dfrac1{\xi}\right).$$
Note that
$$t(\xi) \geq -\dfrac1{\chi} \log \dfrac{\delta_0(\xi)}K \eqno (3.4)$$
for small positive values of $\xi$. The validity of the first two statements of the lemma is clear.

The mapping $h$ is continuous and smooth on the set ${\cal U}\setminus \{0\}$. Let us prove that it is injective. Note that for any $\hat z$, the points $\hat z$ and $z=h(\hat z)$ belong to the same trajectory of system (3.1). Thus, if ${\hat z}_1\neq {\hat z}_2$ and $h({\hat z}_1)= h({\hat z}_2)$, then $|{\hat z}_1|\neq |{\hat z}_2|$. Since the function $t(\xi)$ is strictly monotonous and the Euclidean norm decreases along solutions of system (3.1),
$|h({\hat z}_1)|\neq |h({\hat z}_2)|$, which gives us a contradiction. Hence, $h$ is a homeomorphism.

Any vector $z\in {\mathbb R}^n\setminus \{0\}$ may be represented as $z=r\varphi$, where $r=|z|$ and $\varphi$ is a unit vector. Then
$$Dh(z)=\dfrac{\partial h}{\partial(r,\varphi)}\dfrac{\partial(r,\varphi)}{\partial z}(z)$$
and
$$D_{h,r,\varphi}=\dfrac{\partial h}{\partial(r,\varphi)}=\left(\dfrac{d(\exp(-t(r^2)P)r)}{dr} \varphi, \exp(-t(r^2)P)r (0,E_{n-1})\right).$$
The last $n-1$ columns of the matrix $D_{h,r,\varphi}$ form a basis of the tangent space to the unit sphere at the point $\varphi$. The first column is a nonzero vector orthogonal to this sphere at the same point.

Hence, the matrix $D_{h,r,\varphi}$ is nondegenerate.

The mapping $h^{-1}$ is, therefore, locally Lipschitz continuous and, consequently, globally Lipschitz continuous over compact sets, that do not contain the origin.

Hence there is a positive number $\rho$ such that for all $0<a<b\rho$ there exists a positive number $c(a,b)$ such that
$|h(\hat z_1)-h(\hat z_2)|\geq c(a,b)|\hat z_1-\hat z_2|$ for all $\hat z_{1,2}$ such that $|\hat z_{1,2}|\in [a,b]$. Consequently, the mapping $h^{-1}$ is smooth everywhere except 0.

It follows from (3.3) and (3.4) that if $\hat z$: $|\hat z|<1$, then
$$|\exp (t(\hat z^2)P)\hat z|\leq |\exp (t(\hat z^2)P)|\leq \delta_0(|\hat z|^2),$$
which implies that $|\hat g(\hat z)|\leq |\hat z|^2$. Consequently, $D \hat g (0)=0$.

However, it is still to be proved that $\hat g$ is smooth in a neighborhood of the origin.

Rewrite the conjugacy $h\circ \hat F_2=F_2 \circ h$ in the form
$$\exp (-t((\hat F_2(\hat z))^2)P)\hat F_2(\hat z)=\exp ((-t(\hat z^2)+1)P)\hat z.\eqno (3.5)$$

Formula (3.5) defines a function $\tau (\hat z)$, smooth everywhere except the origin and such that
$$\hat F_2(\hat z)=\exp(\tau(\hat z)P)\hat z.$$
It follows from (3.5) that
$$-t((\exp(\tau(\hat z)P)\hat z)^2)+\tau(\hat z)=-t(\hat z^2)+1.$$
Consequently, $\tau(\hat z)$ is the solution of the equation
$$1-\tau=t(\hat z^2)-t((\exp(\tau P)\hat z)^2)\eqno (3.6)$$

Since $|Cz|>|z|$ for all $z\neq 0$ and the mappings $h$ and $h^{-1}$ preserve the direction along trajectories of System (3.1), $\tau(\hat z)>0$. On the other hand, since the function $t$ decreases, $\tau(z)<1$.
The right-hand side of (3.6) equals to $-t'(\theta)\left((\exp(\tau P)\hat z)^2-\hat z^2\right)$, $\theta\in [\hat z^2,(\exp(\tau(\hat z) P)\hat z)^2]$. Consequently, $\tau(\hat z)\to 0$ as $\hat z\to 0$; otherwise, the right hand side of (3.6) is unbounded.

Rewrite (3.6) in the following form:
$$t(\hat z^2)-1=t((\exp(\tau(\hat z)P)\hat z)^2)-\tau(\hat z).\eqno (3.7)$$
Let a function $\sigma(x,y)$ be defined by the formula $t(x)-y=t(x+\sigma(x,y))$. Note that for any $M>0$ and $k\in {\mathbb N}\bigcup \{0\}$,
$$\dfrac{\partial^k \sigma(x,y)}{\partial x^k}\rightrightarrows 0 \quad \mbox{as } x\to 0, \quad |y|\le M.$$

Then, applying the function $t^{-1}$ to (3.7), we obtain the equality
$$\hat z^2+\sigma(\hat z^2,1)=(\exp(\tau(\hat z)P)\hat z)^2+\sigma((\exp(\tau(\hat z)P)\hat z)^2,\tau(\hat z)).$$

Thus, the difference $(\exp(\tau(\hat z)P)\hat z)^2-\hat z^2$ is a flat function and, consequently, $D\tau(\hat z)\to 0$ and $D{\hat F}_2(\hat z)\to E$ as $\hat z\to 0$. $\blacksquare$

Make a local transformation of variables $\hat x=H(x)=(y,\hat z)=(y,h(z))$ (this transformation can be extended to a global one). The mapping $F$ in the new coordinates takes the form
$$\hat F(\hat x)=H(F(H^{-1}(x)))=(By, {\hat F}_2(\hat z)).$$
The fixed point at the origin is not hyperbolic any more (the unstable manifold becomes center unstable) but it is still strongly conditionally unstable. The cental (former unstable) manifold of the diffeomorphism $\hat F$ is smooth in a neighborhood of the homoclinic point. Due to item 5 of Lemma 3.3, the homoclinic intersection becomes transverse after application of the transformation. Then the statement of the corollary follows directly from that of Theorem 2.4 (where we take $F=G$). $\blacksquare$

\section{Shadowing in a neighborhood of a nonhyperbolic fixed point}

Let $F$ be a homeomorphism of a metric space $(X,\mbox{\rm dist})$.

As usual, we say that a sequence $\{p_k\in X:\;k\in\ensuremath{\mathbb Z}\}$ is a $d$-pseudotrajectory of $F$ if
$$\mbox{\rm dist}(p_{k+1},f(p_k))<d,\quad k\in\ensuremath{\mathbb Z}.$$
We say that a pseudotrajectory $\{p_k:\;k\in\ensuremath{\mathbb Z}\}$ is $\varepsilon$-shadowed by a point $r$ if
$$\mbox{\rm dist\,}(F^k(r),p_k)<\varepsilon, \quad k\in\ensuremath{\mathbb Z}.$$
We say that $f$ has the standard shadowing property if for any $\varepsilon>0$ we can find a $d>0$ such that any $d$-pseudotrajectory of $F$ is
$\varepsilon$-shadowed by some point.

It is well known (see \cite{krpil42}) that to establish the standard shadowing property on a compact phase space it is enough to show that $F$ has the so-called finite shadowing property:

For any $\varepsilon>0$ we can find a $d>0$ (depending on $\varepsilon$ only) such that if $\{p_k:\;0\leq k\leq m\}$ is a finite $d$-pseudotrajectory, then there is a point $r$ such that
$$\mbox{\rm dist\,}(F^k(p),p_k)<\varepsilon, \quad 0\leq k\leq m. \eqno (4.1)$$

Our goal is to give sufficient conditions under which a homeomorphism has the finite shadowing property on a subset of $X$. In our conditions, we use analogs of Lyapunov functions.

Let us formulate our main assumptions.

We assume that there exist two continuous functions
$$W,\,V:X\times X\to\ensuremath{\mathbb R}_+$$
and a compact subset ${\cal N}$ of $X$ on which conditions (C1)-(C9) stated below are satisfied.

We formulate our conditions not directly in terms of the functions $W$ and $V$ but in terms of some geometric objects defined via these functions. Our main reasoning for the choice of this form of conditions is as follows:

(1) Precisely these conditions are used in the proofs below;

(2) it is easy to check conditions of that kind for particular functions $W$ and $V$ (see the example below).

Fix a number $a>0$ and a point $p\in X$ and let
$$\begin{array}{c} P(a,p)=\{q\in X:\;W(q,p)\leq a,\;V(q,p)\leq a\},\\
Q(a,p)=\{q\in P(a,p):\;\;V(q,p)=a\},\\
T(a,p)=\{q\in P(a,p):\;V(q,p)=0\}.\end{array}$$
Let $b>a$; define one more object $R(b,a,p)=\{q:\;a\leq W(q,p)\leq b,\;V(q,p)\leq a\}$.
(Fig.\, 5).

\begin{figure}\begin{center}
\includegraphics*[width=3in]{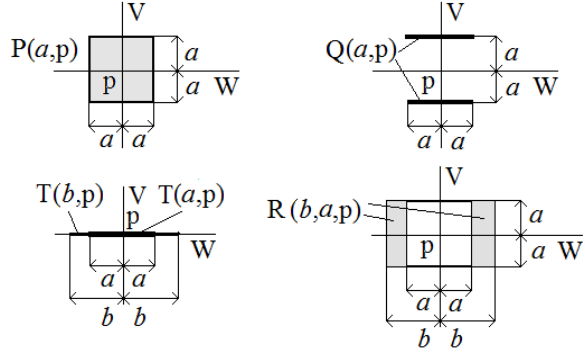}
\end{center}
\caption{\it Sets $P,Q,T$ and $R$.}
\end{figure}

Denote by $B(\varepsilon,p)$ the open $\varepsilon$-ball centered at $p$.

(C1) For any $\varepsilon>0$ there exists a $\delta_0=\delta_0(\varepsilon)>0$ such that $P(\delta,p)\subset B(\varepsilon,p)$ for $p\in{\cal N}$ and $\delta<\delta_0$.

(C2) There exists a $\delta_1>0$ such that if $p\in{\cal N}$ and $\delta<\delta_1$, then there exists a $\Delta>\delta$ and a number $\alpha>0$ such that

(C3) $Q(\delta,p)$ is not a retract of  $P(\delta,p)$;

(C4) $Q(\delta,p)$ is a retract of  $P(\delta,p)\setminus T(\delta,p)$;

(C5) $F(T(\delta,p))\subset \mbox{Int}\,P(\delta,F(p))$;

(C6) $F(P(\delta,p))\subset \mbox{Int}\,P(\Delta,F(p)),\quad F^{-1}(P(\delta,F(p)))\subset \mbox{Int}\,P(\Delta,p)$;

(C7) $F^{-1}(Q(\delta,F(p)))\cap T(\Delta,p)=\emptyset$;

(C8.1) if $q\in P(\Delta,p)$ and $V(q,p)\geq\delta$, then $V(F(q),F(p))>V(q,p)$;

(C8.2) if $q\in P(\Delta,F(p))$ and $W(q,F(p))\geq\delta$, then $$W(F^{-1}(q),p)>W(q,F(p));$$

(C9) there exists a retraction
$$\sigma:\,P(\delta,p)\cup R(\Delta,\delta,p)\to P(\delta,p)$$
such that $V(\sigma(q),p)\geq \alpha V(q,p)$ for $q\in R(\Delta,\delta,p)$.
(Fig.\, 6).

\begin{figure}\begin{center}
\includegraphics*[width=4.5in]{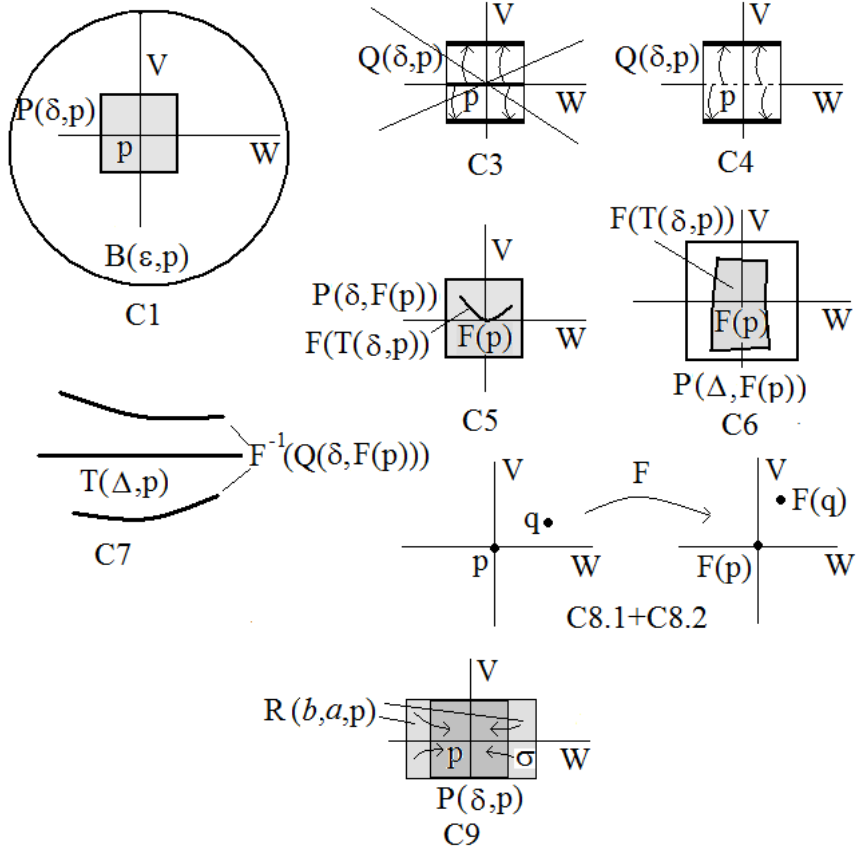}
\end{center}
\caption{\it Conditions (C1)--(C9).}
\end{figure}

\noindent{\bf Theorem 4.1. }{\em Under conditions} (C1)--(C9), $F$ {\em has the finite shadowing property on the set ${\cal N}$: for any $\varepsilon>0$ there exists a $d>0$ such that if $p_0,\dots,p_m$ is a $d$-pseudotrajectory of $F$ belonging to ${\cal N}$, then there is a point $r$ such that inequalities} (4.1) {\em hold.}

In the proof of this statement, we apply the following two lemmas.

Let $p,p'\in{\cal N}$ and $\delta<\delta_1$. We say that condition ${\cal G}(\delta,p,p')$ holds (Fig.\, 7) if
$$F(P)\cap\partial P'\subset Q',\eqno (4.2)$$
$$F(Q)\cap P'=\emptyset, \eqno (4.3)$$
and $Q$ is a retract of the set
$$H=H_1\cup F^{-1}(Q'),$$
where $Q=Q(\delta,p)$, $P=P(\delta,p)$, $Q'=Q(\delta,p')$, $P'=P(\delta,p')$, and $$H_1=P\setminus F^{-1}(\mbox{Int}\,P').$$

\begin{figure}\begin{center}
\includegraphics*[width=3.6in]{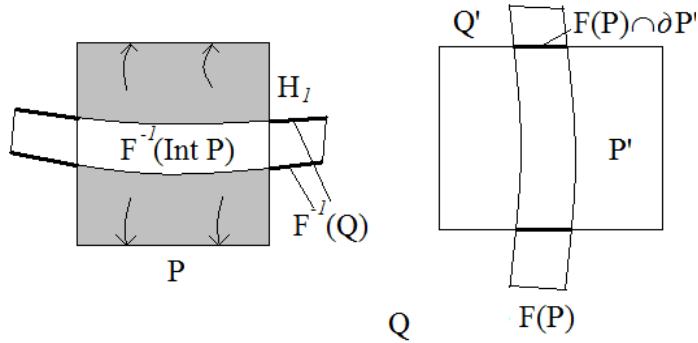}
\end{center}
\caption{\it Condition ${\cal G}(\delta,p,p')$.}
\end{figure}

\noindent{\bf Lemma~4.2. }{\em For any $\delta\in(0,\delta_1)$ there exists a $d_1=d_1(\delta)$ such that if conditions} (C4)-(C9) {\em are satisfied and $\mbox{\rm dist\,}(p',F(p))<d_1$, then condition ${\cal G}(\delta,p,p')$ holds}.

{\bf Proof. } In the proof, we several times select a small $d$ (depending on $\delta$) and then take as $d_1$ the minimum of the selected values of $d$.
Fix a $\delta<\delta_1$.

Condition (C5), the compactness of ${\cal N}$, and the continuity of $F$ imply that there exists a number $c_1=c_1(\delta)<\delta$ such that if $q\in F(T(\delta,p))$, then $W(q,F(p)),V(q,F(p))\leq c_1$. Hence, there exists a $d=d(\delta)$ such that if
$$\mbox{\rm dist\,}(p',F(p))<d,\eqno (4.4)$$
then $W(q,p'),V(q,p')<\delta$, which means that
$$F(T(\delta,p))\subset \mbox{Int}\,P'.\eqno (4.5)$$

A similar reasoning based on Condition (C6) shows that there there exists a $d=d(\delta)$ such that if inequality (4.4) is satisfied, then
$$F(P)\subset P(\Delta,p') \eqno (4.6)$$
and
$$F^{-1}(P')\subset P(\Delta,p). \eqno (4.7)$$

In particular, inclusion (4.7) implies that
$$F^{-1}(Q')\subset P(\Delta,p).\eqno (4.8)$$
Let us show that there exists a $d=d(\delta)$ such that if inequality (4.4) is satisfied, then
$$F^{-1}(Q')\subset P(\delta,p)\cup R(\Delta,\delta,p).\eqno (4.9)$$
Let $S=\{q\in P(\Delta,p):\;V(q,p)\geq \delta\}$. Since the set $S$ is compact, it follows from Condition (C8.1) that there exists a number $c_2=c_2(\delta)>\delta$ such that if $q\in S$, then $V(F(q),F(p))\geq c_2$. Hence, there exists a $d=d(\delta)$ such that if inequality (4.4) is satisfied, then $V(F(q),p')>\delta$ for $q\in S$, which implies that (4.3) is satisfied and $F(S)\cap Q'=\emptyset$. Now inclusion (4.9) follows from inclusion (4.8).

Clearly, Condition (C8.2) (combined with inclusion (4.6)) implies that there exists a $d=d(\delta)$ such that if inequality (4.4) is satisfied, then inclusion (4.2) holds.

Similarly, it follows from Condition (C7) that there exists a number $c_3=c_3(\delta)>0$ such that if $q\in F^{-1}(Q)$, then $V(q,p)\geq 2c_3$. Hence, there exists a $d=d(\delta)$ such that if inequality (4.4) is satisfied, then $V(q,p)\geq c_3$ for $q\in F^{-1}(Q')$.

Apply Condition (C9) to find a retraction $\sigma:\,P(\delta,p)\cup R(\Delta,\delta,p)\to P(\delta,p)$
such that
$$V(\sigma(q),p)\geq \alpha c_3,\quad q\in F^{-1}(Q').\eqno (4.10)$$
The set $U=\sigma(F^{-1}(Q'))$ is compact, and inequality (4.10) implies that
$$U\cap T=\emptyset,\eqno (4.11)$$
where $T=T(\delta,p)$.

By Condition (C4), there exists a retraction $\rho_0$ of $P\setminus T$ to $Q$. Relations (4.5) and (4.11) imply that $H_1\cup U\subset P\setminus T$. Hence, the restriction of $\rho=\rho_0\circ\sigma$ to $H$ is the required retraction $H\to Q$. $\blacksquare$

\noindent{\bf Lemma~4.3. } {\em Let $p_0,\dots,p_m$ be points in ${\cal N}$ such that Condition ${\cal G}(\delta,p_k,p_{k+1})$
holds for $k=0,\dots,m-1$. Then there exists a point $r\in P(\delta,p_0)$ such that $F^k(p)\in P(\delta,p_k)$ for $k=1,\dots,m$.}

{\bf Proof. } Consider the sets (Fig.\, 8)
$$A_k=P(\delta,p_k)\setminus\bigcap_{l>k}F^{-(l-k)}
(\mbox{Int}\,P(\delta,p_l)),\quad k=0,\dots, m-1.$$
It follows from (4.3) that $F(Q(\delta,p_k))\cap P(\delta,p_{k+1})=\emptyset$. Hence, $Q(\delta,p_k)\subset A_k$.

\begin{figure}\begin{center}
\includegraphics*[width=3in]{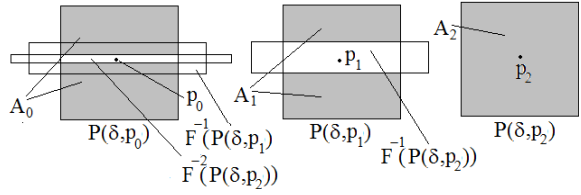}
\end{center}
\caption{\it Sets $A_k$.}
\end{figure}

We claim that there exist retractions $\rho_k:\,A_k\to Q(\delta,p_k),\quad k=0,\dots,m-1$. This is enough to prove our lemma since the existence of $\rho_0$ means that
$$\bigcap_{l=0}^m F^{-l}(\mbox{Int}\,P(\delta,p_l))\neq\emptyset$$
(otherwise there exists a retraction of $P(\delta,p_0)$ to $Q(\delta,p_0)$, which is impossible by Condition (C3)).

The existence of $\rho_{m-1}$ is obvious since condition ${\cal G}(\delta,p_{m-1},p_{m})$ implies the existence of a retraction
$$A_{m-1}\cup F^{-1}(Q(\delta,p_m))\to Q(\delta,p_{m-1}).$$

Let us assume that the existence of retractions $\rho_{k+1},\dots,\rho_{m-1}$ has been proved. Let us prove the existence of $\rho_k$. For brevity, we denote $P_k=P(\delta,p_k), Q_k=Q(\delta,p_k), P_{k+1}=P(\delta,p_{k+1})$, and $Q_{k+1}=Q(\delta,p_{k+1})$.

Note that the definition of the sets $A_k$ implies that
$$A_k\cap F^{-1}(P_{k+1})\subset F^{-1}(A_{k+1}).\eqno (4.12)$$

Define a mapping $\theta$ on $A_k$ by setting $\theta(q)=F^{-1}\circ\rho_{k+1}\circ F(q)$, $q\in A_k\cap F^{-1}(P_{k+1})$,
$\theta(q)=q$, $q\in A_k\setminus F^{-1}(P_{k+1})$. Inclusion (4.12) shows that the mapping $\theta$ is properly defined.

Let us show that this mapping is continuous. Clearly, it is enough to show that $\rho_{k+1}(r)=r$ for $r\in F(A_k\cap F^{-1}(\partial P_{k+1}))$. For this purpose, we note that
$$ F(A_k\cap F^{-1}(\partial P_{k+1}))=F(A_k)\cap\partial P_{k+1} \subset F(P_k)\cap\partial P_{k+1}\subset Q_{k+1}$$
(we refer to inclusion (4.2)) and $\rho_{k+1}(r)=r$ for $r\in Q_{k+1}$.

Clearly, $\theta$ maps $A_k$ to the set
$$[P_k\setminus F^{-1}(P_{k+1})]\cup F^{-1}(Q_{k+1}).\eqno (4.13)$$

Condition ${\cal G}(\delta,p_k,p_{k+1})$ implies that there exists a retraction $\rho$ of (4.13) to $Q_k$. It remains to note that $\theta(q)=q$ for $q\in Q_k$ due to Eq. (4.3). Thus, $\rho_k=\rho\theta:\,A_k\to Q_k$
is the required retraction. The lemma is proved. $\blacksquare$

To complete the proof of the main theorem, we take an arbitrary $\varepsilon>0$, refer to Condition (C1), find the corresponding $\delta_0$, and apply Lemmas~4.2 and 4.3.

\section{Model example}

Let us apply the main result to a diffeomorphism of a neighborhood of the origin of the plane ${\mathbb R}^2$ of the form
$$F(x,y)=(x-x^m+X(x,y),y+y^n+Y(x,y)),$$
where $m>1$ and $n>1$ are odd natural numbers and $X$ and $Y$ are analytic functions with zero Taylor polynomials of degree $m$ and $n$, respectively.
(Of course, it is enough to assume that $X$ and $Y$ are smooth functions whose Taylor series have the corresponding properties; we work with analytic functions to simplify presentation).

Clearly, the origin is a nonhyperbolic fixed point of $F$ (both eigenvalues of the Jacobi matrix $DF(0,0)$ equal 1). Note that the "truncated"\ mapping,
$(x,y)\mapsto (x-x^m,y+y^n)$ is weakly contracting in the $x$ direction and weakly expanding in the $y$ direction.

Let $z,v\in{\mathbb R}$; represent $(z+v)^{2k+1}-z^{2k+1}=vZ_{2k}(z,v)$ for integer $k\geq 0$.

Then $Z_{2k}$ is a positive definite form of degree $2k$ such that
$$Z_{2k}(z,v)\geq (2k+1)z^{2k}. \eqno (5.1)$$

We denote by $p=(p_s,p_u)$ the coordinate representation of a point $p\in{\mathbb R}^2$ and consider the functions
$$W(q,p)=|q_s-p_s|\quad\mbox{and}\quad V(q,p)=|q_u-p_u|.$$

Clearly, condition (C1) is satisfied for any $p\in{\mathbb R}^2$.

Let us show that, under proper smallness conditions on the functions $X$ and $Y$ in a compact neighborhood ${\cal N}$ of the origin, there exists a number $\delta_1$ such that conditions (C2)-(C9) are satisfied with $\alpha=1$.

Take any compact neighborhood ${\cal N}$ of the origin, any $p\in{\cal N}$, and any $\delta>0$; clearly, conditions (C3), (C4), and (C9) with any $\Delta>\delta$ are satisfied for them (take in (C9) the retraction $\sigma$ along the lines $y=\mbox{\rm const}$).

Now we note that, in a fixed neighborhood of the origin,
$$F^{-1}(x,y)=(x+x^m+\Xi(x,y),y-y^n+H(x,y)),$$
where $\Xi$ and $H$ are analytic functions having the same properties as $X$ and $Y$, respectively.

Since $F$ is a diffeomorphism, there exists a neighborhood ${\cal N}_0$ of the origin, a number $K>1$, and a $\delta_1>0$ such that if $p\in{\cal N}_0$ and $\delta\in(0,\delta_1)$, then condition (C6) holds with $\Delta=K\delta$.

Now we show that if ${\cal N}$ is a compact neighborhood of the origin contained in ${\cal N}_0$, proper smallness conditions on $X,Y,\Xi,H$ are satisfied in ${\cal N}$, and $\delta_1$ is small enough, then conditions (C5), (C7), and (C8) hold as well. Let us write down the required conditions of smallness for $X$ and $Y$.

For an arbitrary $\varepsilon>0$, we can find a neighborhood ${\cal N}$ and a number $\delta_1$ such that
$$|X(p_s+v_s,p_u+v_u)-X(p_s,p_u)|\leq\varepsilon(|v_s|+|v_u)|\eqno (5.1)$$
and
$$|Y(p_s+v_s,p_u+v_u)-Y(p_s,p_u)|\leq\varepsilon(|v_s|+|v_u)|\eqno (5.2)$$
for $p\in{\cal N},|v_s|,|v_u|\leq K\delta_1$.

We assume, in addition, that
$$\begin{array}{c}
|X(p_s+v_s,p_u)-X(p_s,p_u)|\leq (m-1)p_s^{m-1}|v_s|,\\ |X(v_s,p_u)-X(0,p_u)|\leq\frac{1}{2}|v_s|^m,
\end{array}
 \eqno (5.3)$$
$$|Y(p_s+v_s,p_u+v_u)-Y(p_s,p_u)|\leq\frac{n-1}{K+1}p^{n-1}_u(|v_u|+|v_s|),\eqno (5.4)$$
and
$$|Y(p_s+v_s,v_u)-Y(p_s,0)|\leq\frac{1}{2}|v_u|^n \eqno (5.5)$$
for $p\in{\cal N},|v_s|,|v_u|\leq K\delta_1$.

Let us show that all the required conditions are satisfied.

Condition (C7). Take $q\in T(\Delta,p)$, let $q=p+(v_s,v_u)$. Then $|v_s|\leq\Delta=K\delta$ and $v_u=0$. Thus,
$$V(F(q),F(p))=|Y(p_s+v_s,p_u)-Y(p_s,p_u)|.$$
It follows from (5.2) that $|Y(p_s+v_s,p_u)-Y(p_s,p_u)|\leq\varepsilon|v_s|\leq\varepsilon K\delta\leq \delta/2$
if $\varepsilon$ is small enough.
Thus,
$$V(F(q),F(p))\leq\frac{1}{2}\delta.\eqno (5.6)$$
This shows that condition (C7) holds.

Condition (C5). Take $q\in T(\delta,p)$; then $q=p+(v_s,v_u)$ with $|v_s|\leq\delta$ and $v_u=0$.

We note that $W(F(q),F(p))=|v_s-v_sZ_{m-1}(p_s,v_s)+X(p_s+v_s,p_u)-X(p_s,p_u)|$.

If $v_s=0$, then $W(F(q),F(p))=0<\delta$.

If $v_s\neq 0$ and $p_s=0$, then it follows from the second condition in (5.3) that
$$W(F(q),F(p))=|v_s-v^m_s+X(v_s,p_u)-X(0,p_u)|<|v_s|\leq\delta.$$

If $v_s\neq 0$ and $p_s\neq 0$, then $Z_{m-1}(p_s,v_s)\geq mp_s^{m-1}$ by (5.1), and we apply the first condition in (5.3) to conclude that $W(F(q),F(p))<|v_s|\leq\delta$.

These inequalities combined with (5.6) show that condition (C5) holds.

Condition (C8.1). Take $q\in P(\Delta,p)$, let $V(q,p)\geq\delta$ and $q=p+(v_s,v_u)$. Then $|v_s|\leq\Delta=K\delta$ and $\delta\leq|v_u|\leq K\delta$.

We note that
$$V(f(q),f(p))=|v_u+v_uZ_{n-1}(p_u,v_u)+Y(p_s+v_s,p_u+v_u) -Y(p_s,p_u)|.$$

If $p_u=0$, then, by condition (5.5),
$$V(f(q),f(p))=|v_u+v^n_u+Y(p_s+v_s,v_u)-Y(p_s,0)|>|v_u|=V(q,p).$$

If $p_u\neq 0$, then we note that $|v_u|Z_{n-1}(p_u,v_u)\geq n|v_u|p^{n-1}_u$ by (5.1). Since $|v_s|\leq K\delta\leq K|v_u|$ and inequality (5.4) holds,
$$\begin{array}{c}|Y(p_s+v_s,p_u+v_u)-Y(p_s,p_u)|\leq\frac{n-1}{K+1}p^{n-1}_u(|v_s|+|v_u|)\leq \\ (n-1)|v_u|p^{n-1}_u,\end{array}$$
and we conclude that $V(F(q),F(p))>|v_u|=V(q,p)$. Thus, condition (C8.1) holds.

Similar assumptions on the representation of $F^{-1}$ allow us to ensure condition (C8.2) (we do not write them explicitly).

\section{Discussion}

There are at least three completely different approaches used in the theory of stability of dynamical systems. The first one, call it "the $C^1$ approach"\ is the most widely used. $C^1$ small perturbations of dynamical systems are considered, and lots of results on structural stability and robustness have been proved. For example, one can mention the Perron theorem, the robustness of hyperbolic invariant sets, or the Pugh closing lemma. For some purposes, for example, $C^1$ linearization, a different approach, which may be called $C^N$, or $C^\infty$, or even $C^{\omega}$ approach is used. As classical example of using these methods one may mention the KAM theory and the theory of normal forms. The main idea of these methods is to use "very smooth"\ perturbations satisfying some additional conditions, like no-resonance ones. However, these methods prove to be very effective in the nonhyperbolic case. The third way may be called "the $C^0$ approach"\ and is used, for example, in Sharkovsky's theory of 1D mappings. The good point of this approach is that no hyperbolicity condition is necessary any more. This method is based on $C^0$ results like the Brauer (or, more general the Schauder) theorem. It is a very good way to prove some existence or shadowing-type results. However, one cannot expect the uniqueness of any structure obtained via this method. The main idea is that a perturbed function, $C^0$ close to the initial one, does have in general the set of zeros, close to ones of the initial function, but one can not expect these zeros are locally unique.

In this paper, we develop a new approach, which is a combination of the first and third ones. In the first part of our paper, we study the chaotic dynamic in a neighborhood of a nonhyperbolic fixed point. We assume that there is a Lyapunov function which provides the total instability of the considered mapping reduced to the center unstable manifold. These assumptions allow us to construct an invariant set of local disks of center unstable manifolds. Applying the Brauer theorem we establish the existence of an infinite set of periodic points of the considered mapping. Our conditions, sufficient for existence of an infinite set of periodic points, are robust, i.e., they can be used in numerical simulations.

These methods, using the Lyapunov functions, have been developed at the second part of the paper. The main idea of this part is that we do not need any hyperbolicity conditions to get the local shadowing property in a neighborhood of a nonhyperbolic fixed point. More precisely, the hyperbolicity condition may be replaced by the existence of two Lyapunov functions, whose level surfaces satisfy some explicitly written topological conditions. This result gives a new glance on possible generalization of hyperbolicity.

In all the cases the obtained structures are not unique.

\section{Conclusion}

A new approach to study the dynamics in neighborhoods of nonhyperbolic periodic points is developed. A new type of chaos in a neighborhood of a nonhyperbolic homoclinic point (Theorem 2.4) and a hyperbolic cubic type homoclinic tangency was introduced (Theorem 3.1). The nonhyperbolic conditions of shadowing in a neighborhood of a fixed point of a homeomorphism were presented (Theorem 4.1). A nonhyperbolic analogue of the $\lambda$~-- lemma was proved (Lemma 2.1).

\bigskip

\textbf{Acknowledgements.} This work was supported by the UK Royal Society, by the Russian Federal Program "Scientific and pedagogical cadres", grant no.  2010-1.1-111-128-033 and by the Chebyshev Laboratory (Department of Mathematics and Mechanics, Saint-Petersburg State University) under the grant 11.G34.31.0026 of the Government of the Russian Federation.

\bigskip

\noindent\textbf{Appendix. Proof of Lemma 2.1.}

\bigskip

Fix a disk $D\in {\cal D}$. Let $W=W^s_{loc}\times W^{cu}_{loc}$ and let $r$ be a point of intersection of $D$ with $W^s_{loc}$. Consider a neighborhood $U=\{(y,z):|y|<\varepsilon^y_2,|z|<\varepsilon^z_2\}$ such that $F(\overline{U})\bigcup \overline{U}\bigcup F^{-1}(\overline{U})\subset W$. Consider local coordinates $x=(y,z)$ in $U$. We may take the neighborhood $U$ such that any admissible disk is a subset of $U$. The proof of Lemma 2.1 can be reduced to that of the following auxiliary statement.

\textbf{Lemma A1.}\emph{There exists a $\delta>0$ such that for any $\varepsilon>0$ there is an integer $m=m(\varepsilon)$ independent on the choice of $D\in {\cal D}$)and a set of smooth disks $\nu_0$, \ldots, $\nu_m$ such that
$$r\in \nu_0\subset D, \qquad r_k:=F^k(r)\in \nu_k\subset F(\nu_{k-1})\bigcap U,\quad k=1,\ldots, m, \eqno (A1)$$
and the disk $\nu_m$ is the graph of a function $y=\beta (z)$, $|z|<\delta$, where
$$|\beta(z)|<\varepsilon \quad\mbox{and}\quad |D\beta(z)|<\varepsilon \quad \mbox{for all}\quad z \quad \mbox{with} \quad |z|<\delta.\eqno (A2)$$}


Let us show how the statement of Lemma 2.1 may be deduced from that of Lemma A1.

The disk $\nu_m$ is the image $b(D_0)$ of the disk
$$D_0=\{x=(0,z)\in W^{cu}_{loc}:|z|\le\delta\}$$
under a mapping $b:(0,z)\to (\beta(z),z)$. It follows from (A1) that $\nu_m\subset F^m(D)$ and, due to (A2), the value $\mbox{\rm dist}_1({\mathop{\rm id}},b)$ can be done arbitrarily small. This proves the statement of Lemma 2.1.

Let us prove Lemma A1. The diffeomorphism $F$ in the neighborhood $U$ can be represented as
$$F(x)=Ax+f(x).$$
On the other hand, with respect to the splitting $x=(y,z)$, we can write $F(x)=(F_1(y,z), F_2(y,z))$, where
$$F_1(y,z)=By+f_1(y,z) \quad \mbox{and} \quad F_2(y,z)=Cy+f_2(y,z).$$

Since the local stable and unstable manifolds coincide with the coordinate spaces in a neighborhood of $0$,
$$\begin{array}{c}
f_2(y,0)=0, \quad \dfrac{\partial f_2}{\partial y}(y,0)=0\qquad \mbox{for all}\quad (y,0)\in W^s_{loc};\\[7pt]
f_1(0,z)=0, \quad \dfrac{\partial f_1}{\partial z}(0,z)=0\qquad \mbox{for all}\quad (0,z)\in W^{cu}_{loc}.
\end{array}\eqno (A3)$$

Fix a $\kappa>0$ and a $\sigma\in (a_0, 1/\|C^{-1}\|)$ such that
$$a=a_0+\kappa<\sigma,\qquad b=b_0-\kappa>\sigma,\qquad \kappa<\dfrac{(b-a)^2}8.\eqno (A4)$$
The neighborhood $U$ may be taken so that
$$|Df(x)|<\kappa \eqno (A5)$$
for all $x\in U$.

Consider a unit vector $v=(v^y,v^z)$ tangent to the disk $\nu_0$ at the point $r$. Since the disk $\nu_0$ is transverse to $W^s_{loc}$ at the point $r$,  $v^z\neq 0$. Consider the so-called inclination $\lambda_0:=|v^y|/|v^z|$ of the vector $v$. Clearly, there exists $\Lambda>0$ such that $\lambda_0<\Lambda$ for any $D\in{\cal D}$ and any nonzero $v\in T_r \nu_0$. Since $r\in W^s_{loc}$,
$$DF(r)=\left(
\begin{array}{cc}
B+\dfrac{\partial f_1}{\partial y}(r) & \dfrac{\partial f_1}{\partial z}(r) \\[5pt]
0 & C+\dfrac{\partial f_2}{\partial z}(r) \\
\end{array}
\right)$$

If $v_1=DF(r)v=(v_1^y,v_1^z)$, then
$$v^y_1=\left(B+\dfrac{\partial f_1}{\partial y}(r)\right)v^y+\dfrac{\partial f_1}{\partial z}(r) v^z \quad \mbox{and} \quad
v^z_1=\left(C+\dfrac{\partial f_2}{\partial z}(r)\right)v^z.$$
The second equality can be rewritten in the form
$$C^{-1}v_1^z=v^z+C^{-1}\dfrac{\partial f_2}{\partial z}(r)v^z.$$
Therefore,
$$|v_1^y|\leq a |v^y|+\kappa |v^z|\quad \mbox{and} \quad |v_1^z|\geq b |v^z|. \eqno (A6)$$
Dividing the first inequality of (A6) by the second one, we see that
$$\lambda_1:=\dfrac{|v_1^y|}{|v_1^z|}\leq \dfrac{a}{b}\lambda_0+\kappa.$$

Denoting the inclinations of the vectors $v_j=DF^j(r)v$ by $\lambda_j$ for $j=1,\ldots,m$, we obtain the estimates
$$\lambda_2\leq \left(\dfrac{a}{b}\right)^2\lambda_0+\dfrac{a\kappa}{b}+\dfrac{\kappa}b, \ldots$$
$$\lambda_m\leq \left(\dfrac{a}{b}\right)^m\lambda_0+\sum_{j=1}^{m-1}\left(\dfrac{a}{b}\right)^j \dfrac{\kappa}b
\leq \left(\dfrac{a}{b}\right)^m\lambda_0+ \dfrac{\kappa}{b-a}.\eqno (A7)$$
Take $m_1\in {\mathbb N}$ so large that
$$\left(\dfrac{a}{b}\right)^{m_1}\Lambda\leq \dfrac{b-a}8.$$
Then it follows from (A4) and (A7) that $\lambda_m\leq (b-a)/4$ for all $v\in T_r \nu_0$, $m\ge m_1$.

We construct the disks $\nu_1=\nu_1(D)$, \ldots, $\nu_{m_1}=\nu_{m_1}(D)$ consecutively choosing
$$\nu_{k+1}= F(\nu_k)\bigcap V_0.$$ Clearly, $r_{k}\in \nu_k$ for all $k=1,\ldots, m_1$.

Since $\nu_{m_1}$ is a smooth disk containing the point $r_{m_1}$ and inequalities (A7) are satisfied, we can choose a $\sigma^0>0$ (independent on the choice of the initial admissible disk $D$) so that
\begin{enumerate}
\item $\mbox{dist}\, (x,W^s_{loc})>\sigma_0$ for any $x\in\partial \nu_{m_1}$;
\item if $\widetilde{\nu}=\nu_{m_1}\bigcup B(\sigma^0,r_{m_1})$, the inequality
$$\lambda<(b-a)/2.\eqno (A8)$$
is satisfied for any $x\in \widetilde{\nu}$ and any $v\in T_x \widetilde{\nu}\setminus \{0\}$.
\end{enumerate}

Define
$$\partial_1=\bigcup_{ D\in {\cal D}} \partial \widetilde{\nu}(D).$$
Clearly, $\partial_1$ is a compact set that does not intersect the local stable manifold. Set $b_1=(b+a)/2$.

Choose $\delta>0$ so that the disk $D_0=\{(0,z):|z|\leq \delta\}$ is a subset of the domain $U$. It follows from (A3) that
$$\dfrac{\partial f_1}{\partial z}(x)=0$$
for all $x\in D$. Let $\varepsilon>0$ be so small that
$$\dfrac{a(b-a)}{2b_1}+\dfrac{\varepsilon (b_1-a)}{2b_1}<\dfrac{b_1-a}2. \eqno (A9)$$
Take a neighborhood $V_0\subset U$ of the disk $D_0$ so small that the inequality
$$\left|\dfrac{\partial f_1}{\partial z}(x)\right|\leq \kappa_1:=\dfrac{\varepsilon(b_1-a)}2$$
is true. Without loss of generality, we may assume that $m_1$ is so large that
$r_{m_1}\in V_0$ and $\widetilde{\nu}\in V_0$.

Consider a point $x\in \widetilde{\nu}$ and a unit vector $v=(v^y,v^z)\in T_x \widetilde{\nu}$. Let $\lambda$ be the inclination of the vector $v$, let
$v_1=(v^y_1,v^z_1)=DF(x)v$, and let $\lambda_1=|v^y_1|/|v^z_1|$.
It follows from the equalities
$$v^y_1=\left(B+\dfrac{\partial f_1}{\partial y}(r)\right)v^y+\dfrac{\partial f_1}{\partial z}(r) v^z \quad\mbox{and}\quad
v^z_1=\dfrac{\partial f_1}{\partial y}(r)v^y+\left(C+\dfrac{\partial f_2}{\partial z}(r)\right)v^z$$
and estimates (A4), (A5), and (A9) that
$|v_1^y|\leq a |v^y|+\kappa |v^z|$ and $|v_1^z|\geq b |v^z|-\kappa |v^y|$.
Consequently,
$$\lambda_1\leq \dfrac{a|v^y|+\kappa_1|v^z|}{b|v^z|-\kappa |v^y|}=\dfrac{a\lambda+\kappa_1}{b-\kappa\lambda}.$$

Since $\kappa\in (0,1)$, it follows from inequality (A8) that $b-\kappa\lambda>b_1$. Inequality (A9) and the definition of the number $\kappa_1$ imply the inequalities
$$\lambda_1\leq \dfrac{a}{b_1}\lambda+\dfrac{\kappa_1}{b_1}<\dfrac{b-a}2.$$

Suppose that a point $x\in \nu_m$ is such that $F^j(x)\in V_0$ for all $j=1,\ldots,m$. If $\lambda_j$ are the inclinations of the vectors $v_j=DF^j(x) v$, $j=1,\ldots, m$, then, iterating the above estimate, we obtain the inequality
$$\lambda_m\leq \left(\dfrac{a}{b_1}\right)^m\lambda+\sum_{j=1}^{m-1}\left(\dfrac{a}{b_1}\right)^j \dfrac{\kappa}{b_1}
\leq \left(\dfrac{a}{b_1}\right)^m\lambda_0+ \dfrac{\kappa_1}{b_1-a}.\eqno (A10)$$

It follows from inequalities (A9) and (A10) that there is a number $m_2$ such that $\lambda_m<\varepsilon$ for all $m\geq m_2$ and $D\in {\cal D}$.

To finish the proof of Lemma A1, we prove that there is an index $m_3>0$ such that the image $\pi_z(F^m(\widetilde{\nu}))$ covers the set $N$ for any $m\geq m_3$.
Here $\pi_z(y,z)=z$.

First we prove that there exist a disk $\Delta=\Delta({\cal D})$ that is covered by every image $\pi_z(F^m(\widetilde{\nu}))$, $m\geq m_3$. If this is not true, there are sequences of boundary points $x_k\in \partial_1$ and of integers $l_k\to\infty$ such that
$$\pi_z(F^{l_k}(x_k))\to 0.\eqno (A11)$$
Since the function $V$ is nonzero over the compact set $\partial_1$, relation (A11) contradicts the strong conditional instability of the fixed point.

Since the fixed point $0$ of the mapping $F|_{W^{cu}_{loc}}$ is a source, one can take a disk $N_0\subset W^{cu}_{loc}$, $0\in N_0$, such that for any disk $\Delta$ with $0\in \Delta\subset W^{cu}_{loc}$ there exists a number $m_4=m_4(\Delta)\in {\mathbb N}$ such that $N_0\subset F^{m_4}(\Delta)$.

The inclinations $|v^y|/|v^z|$ are uniformly small for all $x\in F^m(\widetilde{\nu}(D))$, $D\in {\cal D}$. Consequently, the mappings $h_m:N_0\to F^m(\widetilde{\nu})$, inverse to $\pi_z$, are well-defined and close to identical embeddings provided that $m$ is sufficiently large.

Take $\varrho>0$, $N\subset N_0$, $0\in N$, such that for any embedding $\beta:\Delta\to {\mathbb R}^n$ satisfying the estimate
$$\|\beta-\mbox{\rm id}\|_{C^1}<\varrho, \eqno (A12)$$
$\pi_z(F^{m_4}(\beta(\Delta)))\supset N$ and $\|{\pi_z|_{F^{m_4}(\beta(\Delta))}}^{-1}-\mbox{\rm id}\|_{C^1}<\varepsilon$ in the space of embeddings of the disk $$\pi_z(F^{m_4}(\beta(\Delta)))$$ to ${\mathbb R}^n$.

Then we select an integer $m_5$ such that the disks $F^m(\widetilde{\nu})$, $m>m_5$ contain the images $\beta_m(\Delta(D))$ of mappings $\beta_m$, satisfying (A12). To finish the proof, it suffices to take $h_m={{\pi_z|_{F^{m_4}(\beta_m(\Delta))}}^{-1}}|_N$. $\blacksquare$

\end{document}